\documentclass[10pt]{amsart}
\usepackage{mathptmx}
\usepackage{amsmath}
\usepackage{amssymb}
\usepackage{array}
\usepackage{geometry}
\usepackage[bookmarks=true,colorlinks=true, pdfstartview=FitV, linkcolor=black, citecolor=blue, urlcolor=black]{hyperref}

\usepackage{color}
\definecolor{DarkRed}{rgb}{0.55,.00,0.2}
\definecolor{DarkGrey}{rgb}{0.35,.35,0.35}

\theoremstyle{definition}

\theoremstyle{remark}

\numberwithin{equation}{section}

\begin{document}

\title{New summation and transformation formulas of the  Poisson,  M\"{u}ntz, M\"{o}bius   and Voronoi type }

\author{Semyon Yakubovich$^\ast$}
\thanks{$^\ast$ E-mail: syakubov@fc.up.pt}
\maketitle

\markboth{\rm \centerline{ Semyon   Yakubovich}}{}
\markright{\rm \centerline{New summation and transformation formulas  }}
\begin{center}{\it Department of Mathematics,  Faculty of Sciences,\\   University of Porto,  Campo Alegre st., 687,  4169-007 Porto,    Portugal}\end{center}

\begin{abstract} {\noindent Basing on  properties of the Mellin transform and  Ramanujan's identities, which represent  a ratio of products of Riemann's zeta- functions of different arguments in terms of the  Dirichlet  series of  arithmetic functions, we obtain  a number of the Poisson,   M\"{u}ntz, M\"{o}bius  and Voronoi type summation formulas.  The corresponding analogs of the M\"{u}ntz operators are investigated.  Interesting and curious particular cases of summation formulas involving arithmetic functions  are exhibited. Necessary and sufficient conditions for the validity of the Riemann hypothesis are derived.}
\end{abstract}

\vspace{6pt}


\vspace{6pt}


\section {Introduction. Classical summation formulas}

Let $f: \mathbb{R}_+ \to \mathbb{C}$ be a complex-valued function  whose Mellin's transform is defined by the integral \cite{tit}
$$f^*(s)= \int_0^\infty  f(x) x^{s-1} dx,\quad   s= \sigma + it,\  \sigma, \ t \in \mathbb{R},\eqno(1.1)$$
 which is well defined under certain conditions given below.   The inverse Mellin transform is defined accordingly
 $$f(x)=  {1\over 2\pi i} \int_{\sigma-i\infty}^{\sigma+i\infty} f^*(s) x^{-s} ds,\quad x >0\eqno(1.2)$$
 as well as the Parseval  equality   for two functions $f, g$ and their Mellin transforms  $f^*,\ g^*$

 $$\int_0^\infty  f(xy) g(x) dx =  {1\over 2\pi i} \int_{\sigma-i\infty}^{\sigma+i\infty} f^*(s) g^* (1-s) y^{-s} ds,\  y >0.\eqno(1.3)$$
 Our results will be based on the properties and series representations of  the familiar  Riemann zeta-function \cite{titrim} $\zeta(s)$, which satisfies  the  functional equation
$$\zeta(s)= 2^s \pi^{s-1} \sin\left({\pi s\over 2}\right)\Gamma(1-s)\zeta(1-s),\quad s=\sigma +it,\eqno(1.4)$$
where $\Gamma(z)$ is Euler's gamma function.  In the right half-plane  ${\rm Re}  s  > 1$ it is represented by the absolutely and uniformly convergent series with respect to $t \in \mathbb{R}$
$$\zeta(s) =\sum_{n=1}^\infty \frac{1}{n^{s}},\eqno(1.5)$$
and by the uniformly convergent series
$$(1- 2^{1-s})\zeta(s) =\sum_{n=1}^\infty \frac{(-1)^{n-1}}{n^{s}},\   {\rm Re} s > 0.\eqno(1.6)$$
Moreover,  transformation and summation formulas, which will be derived in the sequel are generated by the  Ramanujan identities (see in \cite{titrim}, \  \cite{kar} ) involving  arithmetic functions and ratios of Riemann's zeta-functions of different arguments, namely
$$\frac{\zeta^2(s)}
{\zeta(2s)} =\sum_{n=1}^\infty \frac{ 2^{\omega(n)}}{n^{s}}, \quad
{\rm Re} s > 1,\eqno(1.7)$$
$$\zeta^2(s) =\sum_{n=1}^\infty \frac{ d(n)}{n^{s}}, \quad
{\rm Re} s > 1,\eqno(1.8)$$
$$\zeta^k(s) =\sum_{n=1}^\infty \frac{ d_k(n)}{n^{s}}, \quad
{\rm Re} s > 1,\  k=2,3,4,\dots,\eqno(1.9)$$
$$\frac{1}{\zeta(s)} =\sum_{n=1}^\infty \frac{\mu(n)}{n^{s}}, \quad {\rm Re} s
> 1,\eqno(1.10)$$
$$\frac{\zeta(s)}
{\zeta(2s)} =\sum_{n=1}^\infty \frac{\left|\mu(n)\right|}{n^{s}},
\quad {\rm Re} s > 1,\eqno(1.11)$$
$$\frac{\zeta(2s)}
{\zeta(s)} =\sum_{n=1}^\infty \frac{\lambda(n)}{n^{s}}, \quad {\rm
Re} s > 1,\eqno(1.12)$$
$$\frac{\zeta^3(s)}
{\zeta(2s)} =\sum_{n=1}^\infty \frac{d\left(n^2\right)}{n^{s}},
\quad {\rm Re} s > 1,\eqno(1.13)$$
$$\frac{\zeta^4(s)}
{\zeta(2s)} =\sum_{n=1}^\infty \frac{d^2(n)}{n^{s}}, \quad {\rm Re}
s > 1,\eqno(1.14)$$
$$\frac{\zeta(s-1)}
{\zeta(s)} =\sum_{n=1}^\infty \frac{\varphi(n)}{n^{s}}, \quad {\rm
Re} s > 2,\eqno(1.15)$$
$$\frac{1- 2^{1-s}}{1-2^{-s}}\zeta(s-1) =\sum_{n=1}^\infty \frac{a(n)}{n^{s}},
\quad {\rm Re} s > 2.\eqno(1.16)$$
The arithmetic function  $a(n)$ in (1.16) denotes the greatest odd divisor of $n$,
$d(n)$ in (1.8)  is the Dirichlet divisor function, i.e. the number of divisors of $n$,
including $1$ and $n$ itself.    A more general  function $d_k(n),\ k=2,3,4\dots$ denotes the number of ways of expressing $n$ as a product of $k$ factors and expressions with the same factors in a different order being counted as different.    The divisor function has  the estimate \cite{titrim}  $d(n)=  O(n^\varepsilon), \ n \to \infty, \ \varepsilon
> 0$. The M\"{o}bius function in (1.10), (1.11) is denoted by $\mu(n)$ and $\mu(1)=1,\  \mu(n)= (-1)^k$ if $n$ is the product of $k$ different primes, and $\mu(n)=0$ if $n$ contains any factor to a power higher than the first. The symbol $\omega(n)$ in (1.7) represents the  number of distinct prime factors of $n$ and it behaves as $\omega(n)= O(\log\log n), \ n \to \infty$. By  $\varphi(n)$ in (1.15) the Euler totient function is denoted and  its asymptotic  behavior satisfies $\varphi(n)= O\left(n [\log\log n]^{-1}\right), \ n \to \infty$. Finally, $\lambda(n)$ in (1.12) is  the Liouville function and it has the estimate $|\lambda(n)| \le 1$.

Recently \cite{yakrim}, the author investigated  invertibility of the transformations with arithmetic functions in a  special class, related to the inverse Mellin transform (1.2) (cf. \cite{yal}). In particular,  it concerns the classical M\"{o}bius expansion (see, for instance, in \cite{edw},  Chapter 10,\  \cite{yakrim} )
$$f(x)= \sum_{n=1}^\infty \mu(n)\sum_{m=1}^\infty  f(xnm),\eqno(1.17)$$
which generates the M\"{o}bius transformation
$$\left(\Theta   f\right) (x)=   \sum_{n=1}^\infty  f(xn),\quad x >0.\eqno(1.18)$$

The classical M\"{u}ntz formula 
$$\zeta(s) f^*(s)=  \int_0^\infty x^{s-1} \left[ \left(\Theta f\right) (x) -  {1\over x} \int_0^\infty f(y) dy \right] dx,\eqno(1.19)$$
where $s=\sigma+ it,\ 0 < \sigma < 1$ is proved, for instance,  in \cite{titrim}, Chapter 2 under conditions $f \in C^{(1)}[0,\infty)$, having the asymptotic behavior at infinity $f(x)= O( x^{-\alpha}),\quad \alpha > 1,\ x \to \infty.$ The expression in square brackets in the right-hand side of (1.19) is called the M\"untz operator (cf. \cite{baez}, \  \cite{baez1}, \cite{ivic},\ \cite{alb})
$$(Pf)(x)=   \left(\Theta f\right) (x) -  {1\over x} \int_0^\infty f(y) dy.\eqno(1.20)$$ 

{\bf Definition.} {\it  A function $f(x),\ x \ge 0$ belongs to the M\"{u}ntz type class $\mathcal{M}_\alpha$,  if $f \in C^{(2)}[0,\infty)$ and  its successive derivatives have  the asymptotic behavior at infinity $f^{(j)}(x)= O( x^{-\alpha-j}),\quad \alpha > 1,\ j=0,1,2,\  x \to \infty.$}

The following lemma will give the representation of the M\"{u}ntz operator in terms of the inverse Mellin transform (1.2).

{\bf Lemma 1}. {\it   Let $f \in \mathcal{M}_\alpha$.  Then its Mellin transform $(1.1)$ $f^*(s),\ s= \sigma +it$  is analytic in the strip $0< \sigma < \alpha$ and  belongs  to $L_1(\sigma-i\infty,\ \sigma +i\infty)$ over any vertical line in the strip.  Finally, the M\"{u}ntz operator has the representation
$$(Pf)(x)=   {1\over 2\pi i} \int_{\sigma-i\infty}^{\sigma+i\infty} \zeta(s) f^*(s) x^{-s} ds,\quad x >0,\eqno(1.21)$$
valid for $0 < \sigma < 1$.}

\begin{proof} Indeed, integrating twice by parts in (1.1) and eliminating the integrated terms, we derive
$$f^*(s)=  {1\over s(s+1)}  \int_0^\infty x^{s+1}  f^{(2)}(x) dx,\quad  0 < \sigma <\alpha.\eqno(1.22)$$  
Hence
$$|f^*(s)| \le  {1\over |s(s+1)|}  \left[ \int_0^1 x^{\sigma+1} | f^{(2)}(x)| dx +  \int_1^\infty  O( x^{\sigma-\alpha-1} ) dx \right]
= O(|s|^{-2}),$$ 
which means the analyticity of the Mellin transform (1.1) in the strip $0< \sigma < \alpha$ and its integrability over any vertical line in the strip.  But considering for now $\alpha > \sigma > 1$ and using (1.2), (1.5), we easily get via the change of the order of integration and summation that the M\"{o}bius transformation (1.18) can be written as
$$\left(\Theta   f\right) (x) =   {1\over 2\pi i} \int_{\sigma-i\infty}^{\sigma+i\infty} \zeta(s) f^*(s) x^{-s} ds.\eqno(1.23)$$
This is indeed possible due to the absolute convergence of the iterated series and integral for each $x >0$
$$\int_{\sigma-i\infty}^{\sigma+i\infty} \left(\sum_{n=1}^\infty \frac{1}{n^{\sigma}} \right) | f^*(s) x^{-s} ds| \le
x^{-\sigma} \zeta(\sigma) \int_{\sigma-i\infty}^{\sigma+i\infty} | f^*(s)  ds| < \infty.$$
Hence (see in \cite{tit}), since   $f(x) x^{\sigma-1} \in L_{1}(\mathbb{R}_{+})$,  $f$   can be represented by the absolutely convergent integral (1.2).  Thus, changing the order of integration and summation, we come up with (1.23).

In the meantime, $\zeta(s)$ is bounded for $\sigma > 1$ and for $0 < \sigma \le  1$ it has the behavior (cf. \cite{titrim}, Chapter 3, \cite{kar}, Chapter 2)
$$\zeta(\sigma +it) = O(|t|^{\varepsilon+ (1-\sigma)/2} ),\    \zeta(1  +it) = O(\log t  ),\   |t| \to \infty\eqno(1.24)$$ 
for every positive $\varepsilon$. Moreover, the product $\zeta(s) f^*(s) x^{-s}$ is analytic for each $x >0$ in the strip $0< \sigma < \alpha$ except, possibly,  for a simple pole at $s=1$ with residue  $ f^*(1) x^{-1}$.   The following asymptotic behavior  (see (1.22), (1.24))
$$\zeta(s) f^*(s)  = O(|t|^{\varepsilon - (3+\sigma)/2}),  \quad   0< \sigma < 1,\    |t| \to \infty,$$
$$ \zeta(s) f^*(s)  = O(|t|^{-2} \log t),  \quad    \sigma =1,\    |t| \to \infty,$$
$$\zeta(s) f^*(s)  = O(|t|^{-2}),  \quad   1 < \sigma < \alpha,\    |t| \to \infty,$$
guarantees the absolute integrability in (1.23) over any vertical  contour $\sigma +it,\   |t| \ge t_0 > 0$, where $\sigma$ is lying in the interval $(0,\alpha)$, Therefore via the residue  theorem   it  becomes 
$$\left(\Theta   f\right) (x) =   {1\over 2\pi i} \int_{\sigma-i\infty}^{\sigma+i\infty} \zeta(s) f^*(s) x^{-s} ds +  {1\over x} \int_0^\infty f(y) dy,\ 0 < \sigma < 1.\eqno(1.25)$$
Hence, recalling (1.20),  we get  (1.21), completing the proof of the lemma. 

\end{proof}

{\bf Corollary 1}. {\it  Let $f \in \mathcal{M}_\alpha$ and $f^*(1)= 0$. Then the M\"{o}bius operator $(1.18)$ has the representation
$$\left(\Theta   f\right) (x) =   {1\over 2\pi i} \int_{\sigma-i\infty}^{\sigma+i\infty} \zeta(s) f^*(s) x^{-s} ds, \ x >0,\eqno(1.26)$$
valid for $ 0 < \sigma < 1$.}

\begin{proof}  The proof is immediate from (1.25) and the definition of the Mellin transform (1.1).

\end{proof}

Further, let us consider the familiar Poisson formula \cite{tit}
$$\sqrt x \left[ {1\over 2} (F_c f) (0) + \sum_{n=1}^\infty  (F_c f)(nx)\right] = \sqrt{2\pi\over x} \left[   {1\over 2}  f (0) + \sum_{n=1}^\infty   f \left({2\pi n\over x}\right)\right],\ x >0,\eqno(1.27)$$
where $(F_cf)(x)$ denotes the operator of the Fourier cosine transform
$$(F_cf)(x)=  \sqrt{2\over \pi } \int_0^\infty f(t) \cos (xt) dt.\eqno(1.28)$$

We will give a rigorous proof of the Poisson formula (1.27) in the class $\mathcal{M}_\alpha$,  justifying  the formal method proposed by Titchmarsh in \cite{tit}, Section 2.9.    Precisely, it has 

{\bf Theorem 1}. {\it Let $f \in \mathcal{M}_\alpha$. Then for all $x >0$ the Poisson formula $(1.27)$ holds.}

\begin{proof}  In fact,  since $f\in \mathcal{M}_\alpha$,  integral (1.28) converges absolutely and uniformly on $\mathbb{R}_+$.    Now integrating twice by parts in the integral (1.28), we find
$$(F_cf)(x)=  \sqrt{2\over \pi } {1\over x^2} \left[ f^{(1)}(0) -  \int_0^\infty f^{(2)} (t) \cos (xt) dt \right]= O\left({1\over x^2} \right),\  x \to \infty.\eqno(1.29)$$
Hence,  differentiating with respect to $x$, we obtain
$$(F_cf)^{(1)} (x)=  - 2 \  \sqrt{2\over \pi } {1\over x^3} \left[ f^{(1)}(0) -  \int_0^\infty f^{(2)} (t) \cos (xt) dt \right]$$
$$+ \sqrt{2\over \pi } {1\over x^2}  \int_0^\infty t f^{(2)} (t) \sin (xt) dt = O \left({1\over x^2} \right),\  x \to \infty,$$
where the differentiation is allowed under the integral sign via the absolute and uniform convergence, since $t f^{(2)} (t)= O(t^{-\alpha-1}), \   t \to \infty.$  Similarly, we get  
$$(F_cf)^{(2)} (x)=   O \left({1\over x^2} \right),\  x \to \infty$$
and certainly $(F_cf)(x) \in C^2 (\mathbb{R}_+)$.    Further, in the class $\mathcal{M}_\alpha$ integral (1.28) can be written  in the form
$$(F_cf)(x)=  \sqrt{2\over \pi } {d\over dx} \int_0^\infty f(t) {\sin (xt) \over t} dt.$$
Calling the well-known integral \cite{tit}
$$\int_0^\infty  {\sin t \over t} \ t^{s-1} dt = \cos\left({\pi s\over 2}\right) {\Gamma(s)\over 1-s},\quad  0 < \sigma < 1,$$
and observing that $f(x) \in L_2(\mathbb{R}_+)$ and $\sin x /x $ is square integrable as well, the Parseval equality (1.3) holds and we find the representation 
$$(F_cf)(x)=  \sqrt{2\over \pi } {d\over dx}  {1\over 2\pi i} \int_{\sigma-i\infty}^{\sigma+i\infty}  f^*(1- s) 
  \cos\left({\pi s\over 2}\right) {\Gamma(s)\over 1-s} x^{1-s} ds,\quad 0 < \sigma < 1. $$
 Meanwhile,  the differentiation with respect to $x$ is possible under the integral sign in the latter equality,  because   owing to the Stirling asymptotic formula for the gamma function \cite{erd}, Vol. I
 $$ \cos\left({\pi s\over 2}\right) \Gamma(s) = O( |s|^{\sigma- 1/2}),\   |s| \to \infty.$$
Consequently (see the proof of Lemma 1),  the function
$$f^*(1- s)   \cos\left({\pi s\over 2}\right) \Gamma(s)  = O( |s|^{\sigma- 5/2}),\   |s| \to \infty,\ 0 < \sigma < 1$$
belongs to $L_1(\sigma-i\infty,  \sigma+i\infty)$ and we derive  the representation for all $x >0$
$$(F_cf)(x)=  \sqrt{2\over \pi }  {1\over 2\pi i} \int_{\sigma-i\infty}^{\sigma+i\infty}  f^*(1- s) 
  \cos\left({\pi s\over 2}\right) \Gamma(s) \  x^{-s} ds,\quad 0 < \sigma < 1, \eqno(1.30)$$
However, the Mellin transform (1.1) $f^*(1-s)$ can be written as 
$$f^*(1-s) = {f(0)\over 1-s} + \int_0^1 [ f(t) - f(0)] t^{-s} dt +  \int_1^\infty f(t) t^{-s} dt = {f(0)\over 1-s}$$
$$ +   {1\over 1- s} \int_0^1  f^{(1)} (t) (1- t^{1-s}) dt -   {f(1)\over 1-s} - {1\over 1-s} \int_1^\infty f^{(1)} (t) t^{1-s} dt
= O(|s|^{-1}),\ |s| \to \infty,  \  1 < \sigma < \hbox{min}(\alpha, 2). $$
So, it possibly has a simple pole at the point $s=1$ with the residue $f(0)$. The function  $f^*(1-s)$ can be continued analytically into the strip $1 < \sigma < \hbox{min} (\alpha,  3/2)$,  and as we see it behaves as $O(|s|^{-1})$ at infinity.  Therefore,  writing (1.30) as 
$$(F_cf)(x)=  \sqrt{2\over \pi }  {1\over 2\pi i}  {d\over dx} \int_{\sigma-i\infty}^{\sigma+i\infty}  f^*(1- s) 
  \cos\left({\pi s\over 2}\right) \Gamma(s) \  {x^{1-s} \over 1-s}ds,\eqno(1.31)$$
one finds that the integrand in (1.31) has a simple pole  at the point $s=1$ with the residue $\pi f(0)/2$. Thus moving the contour to the right via the Cauchy theorem,  we deduce
$$(F_cf)(x)=  \sqrt{2\over \pi }  {1\over 2\pi i}  {d\over dx} \int_{\sigma-i\infty}^{\sigma+i\infty}  f^*(1- s) 
  \cos\left({\pi s\over 2}\right) \Gamma(s) \  {x^{1-s} \over 1-s}ds -  \sqrt{\pi\over 2 }    {d\over dx} [f(0)] $$
$$= \sqrt{2\over \pi }  {1\over 2\pi i}  {d\over dx} \int_{\sigma-i\infty}^{\sigma+i\infty}  f^*(1- s) 
  \cos\left({\pi s\over 2}\right) \Gamma(s) \  {x^{1-s} \over 1-s}ds,\   1 < \sigma < \hbox{min}(\alpha,\  3/2).$$
Hence, equality (1.31) is valid  for all $0 < \sigma < \hbox{min}(\alpha, \  3/2)$ except $ \sigma \neq 1$ in the case $f(0)\neq 0$.

Now we are ready to prove the Poisson formula (1.27).  Indeed, recalling the functional equation (1.4) for
the Riemann zeta-function, we substitute its right-hand side  into (1.21) to obtain after a simple change of variables and differentiation under the integral sign 
$$\sum_{n=1}^\infty   f \left(n x\right)  -  {1\over x} \int_0^\infty f(y) dy =   {1\over 2\pi i} {d\over dx}  \int_{1-\sigma-i\infty}^{1-\sigma+i\infty} 2^{1-s} \pi^{-s}  \zeta(s) f^*(1- s) \cos\left({\pi s\over 2}\right)\Gamma(s) {x^{s}\over s} ds.$$
Hence,  shifting the contour to the right in the right- hand  side of the latter equality due to the Cauchy theorem, we encounter a simple pole of the integrand at the point $s=1$ with the residue $f(0)x /2$.  This is permitted,  since the integrand behaves at infinity (see (1.24) and use again the asymptotic Stirling formula for the gamma-function) as 
$$\zeta(s) f^*(1- s) \cos\left({\pi s\over 2}\right){\Gamma(s)\over s} =  O(|t |^{\varepsilon + \sigma/2 - 3}),\ 0 < \sigma < 1,\ |t| \to \infty$$
$$\zeta(s) f^*(1- s) \cos\left({\pi s\over 2}\right){\Gamma(s)\over s}  = O( |t|^{-3/2}\log t),\ \sigma=1, \ |t| \to \infty,$$ 
$$\zeta(s) f^*(1- s) \cos\left({\pi s\over 2}\right){\Gamma(s)\over s}  = O( |t|^{\sigma -5/2}),\ 
1< \sigma < \hbox{min}(\alpha,\    3/2), \ |t| \to \infty.$$ 
Thus integrating over a vertical line with some $\sigma \in (1,  \ \hbox{min}(\alpha,\    3/2) )$, we take into account (1.31), (1.28),  (1.5) and  after the change of  the order of  integration and summation via the absolute convergence, we derive 
$$\sum_{n=1}^\infty   f \left(n x\right)  + {1\over 2} f(0) =   \sqrt{2\pi}  \left[ {1\over 2x} (F_cf)(0) +   {d\over dx } \sum_{n=1}^\infty \int_0^x   (F_c f)\left({2\pi n\over y}\right) {dy \over y} \right],  \  x >0.$$
But the differentiation of the series is allowed owing to the estimate (1.29), which easily shows its convergence for all $ x >0$ and the uniform convergence by $x \in (0,  x_0],\   x_0 >0$ of the series of derivatives.  Therefore, differentiating the series term by term, we come up  with (1.27) subject to  a simple change of variables. 

\end{proof}

Our goal now will be to derive the form of the Voronoi operator similarly to (1.20), (1.21) and to prove the Voronoi summation formula \cite{titrim}.  Concerning  the recent results of the author on this  subject in  $L_2$ see in \cite{yakvor}, \cite{yakind}, \cite{yakkosh}. 

Let us consider the arithmetic transform, involving the divisor function
$$(Df)(x)= \sum_{n=1}^\infty  d(n) f \left(n x\right),\ x >0.\eqno(1.32)$$ 
Assuming that $f \in \mathcal{M}_\alpha$, one can take the Mellin transform (1.1) from both sides of the equality (1.32), where $1 < \sigma < \alpha$. Then changing the order of integration and summation due to the absolute convergence and using  (1.8), we find
$$(Df)^*(s)= \zeta^2(s) f^*(s),\quad  1 < \sigma < \alpha.\eqno(1.33)$$
However, the right hand side of the latter equality belongs to $L_1(\sigma-i\infty,\ \sigma +i\infty)$ because $\zeta(s)$ is bounded and $f^*(s)= O(|s|^{-2}),\ |s| \to \infty.$  Therefore  for all $x >0$ (see (1.2))
$$(Df) (x)= {1\over 2\pi i}  \int_{\sigma-i\infty}^{\sigma+i\infty}  \zeta^2(s) f^*(s) x^{-s} ds.\eqno(1.34)$$
On the other hand, the integrand in (1.34) is analytic in the strip $0 < \sigma< \alpha$ except $s=1$, where it has a double pole.  Moreover, recalling again (1.22), (1.24), we get 
$$\zeta^2(s) f^*(s) = O(|t|^{\varepsilon -\sigma -1}),\ |t| \to \infty, \   0 < \sigma < 1,$$
$$\zeta^2(s) f^*(s) = O(|t|^{-2}\log^2 t),\ |t| \to \infty, \   \sigma=1.$$
Consequently, the Cauchy theorem allows us to shift the contour to the left, taking into account the residue at the double pole $s=1$. It can be calculated, in turn, employing the Laurent series for zeta-function in the neighborhood of $s=1$ \cite{erd}, Vol. I. Therefore, after straightforward calculations we obtain
$${\rm Res}_{s=1} [ \zeta^2(s) f^*(s) x^{-s}] = \int_0^\infty f(xy) (\log y + 2\gamma) dy,$$
where $\gamma$ is the Euler constant.  Hence we arrived at the equality
$$ {1\over 2\pi i}  \int_{\sigma-i\infty}^{\sigma+i\infty}  \zeta^2(s) f^*(s) x^{-s} ds = (Df)(x)-   \int_0^\infty f(xy) (\log y + 2\gamma) dy,\ x >0,\eqno(1.35)$$
which is valid for $0 < \sigma <1$.   Thus we proved

{\bf Lemma 2}. {\it Let  $f \in \mathcal{M}_\alpha$.  Then the Voronoi operator
$$(Vf)(x)=   \sum_{n=1}^\infty  d(n) f \left(n x\right) -   \int_0^\infty f(xy) (\log y + 2\gamma) dy,\ x >0\eqno(1.36)$$
is well defined and represented by the equality $(1.35)$.}   

An interesting corollary follows immediately, taking into account the M\"{u}ntz formula (1.19), the definition of the M\"{u}ntz operator (1.20), its representation (1.21) and the previous lemma.  We have

{\bf Corollary 2}. {\it Let $f(x)$ and its  M\"{u}ntz operator $(Pf)(x)$ belong to $\mathcal{M}_\alpha$. Then the Voronoi operator $(Vf)(x)$ is equal to the iterated M\"{u}ntz operator, i.e.}
$$(Vf)(x)= (P^2f)(x),\quad x >0.\eqno(1.37)$$

\begin{proof} The proof is immediate from the equality (1.35), where the left hand-side is equal to $(P^2f)(x)$ via the M\"{u}ntz formula (1.19) and representation (1.21).

\end{proof}

{\bf Theorem 2}. {\it Let $f \in \mathcal{M}_\alpha$. Then the M\"{u}ntz type formula for the Voronoi operator 
$$\zeta^2(s) f^*(s) = \int_0^\infty  x^{s-1} \left[  \sum_{n=1}^\infty  d(n) f \left(n x\right) -   \int_0^\infty f(xy) (\log y + 2\gamma) dy\right] dx\eqno(1.38)$$
is valid for $0  < \sigma < 1$.} 

\begin{proof}   We have
$$\int_0^\infty f(xy) (\log y + 2\gamma) dy = {2\gamma \over x} \int_0^\infty f(y) dy +  {1 \over x} \int_0^\infty f(y) \log y \  dy-
{\log x  \over x} \int_0^\infty f(y)  dy =  {c_1\over x} - {c_2 \log x\over x},$$
where $c_j,\ j=1,2$ are constants
$$c_1=  \int_0^\infty f(y) (2\gamma+ \log y ) dy,  \quad  c_2=  \int_0^\infty f(y)  dy.$$
Hence, for ${\rm Re}\ s > 1$
$$ \int_0^\infty  x^{s-1}  \sum_{n=1}^\infty  d(n) f \left(n x\right)   \ dx =  \int_0^1  x^{s-1} \left[ \sum_{n=1}^\infty  d(n) f \left(n x\right)  -  {c_1\over x} + {c_2 \log x\over x} \right] \ dx +  {c_1\over s-1} +  {c_2\over (s-1)^2} $$
$$+ \int_1^\infty  x^{s-1}  \sum_{n=1}^\infty  d(n) f \left(n x\right)   \ dx.\eqno(1.39)$$
However, appealing to (1.35),  and moving the contour to the left in its  left-hand side,  we find a $\delta \in (0,1),  \delta < \sigma $ to  establish the following estimate 
$$\left| \sum_{n=1}^\infty  d(n) f \left(n x\right) -   \int_0^\infty f(xy) (\log y + 2\gamma) dy \right| =
\left|  {1\over 2\pi i}  \int_{\delta-i\infty}^{\delta +i\infty}  \zeta^2(s) f^*(s) x^{-s} ds \right| $$
$$\le {x^{-\delta} \over 2\pi }   \int_{\delta-i\infty}^{\delta +i\infty}  |\zeta^2(s) f^*(s)  ds | = C\   x^{-\delta},$$
where $C >0$ is a constant.  Hence $(Vf)(x)= O (x^{-\delta} ),\ x \to 0 $ and the right-hand side of (1.39) is analytic for $\delta < \sigma < \alpha$ (except at $s=1$).  Moreover, when $\sigma  < 1$
$$  {c_1\over s-1}= - c_1 \int_1^\infty x^{s-2} dx,$$
$${c_2\over (s-1)^2} =  c_2 \int_1^\infty x^{s-2} \log x\  dx.$$
Substituting these values into  (1.39), we  take  in mind equality (1.33),  and since $ \delta \in (0,1) $ is arbitrary, come up with (1.38), completing the proof of Theorem 2.  
\end{proof}

Finally in this section we prove by the same method the Voronoi summation formula. We note that recently another alternative proof of this formula was given in \cite{vor}. 

{\bf Theorem 3.}  {\it Let $f(x)$ and $x^{-1} (F_cf)(x^{-1})$ belong to $\mathcal{M}_\alpha$ with $\alpha > 2$. Then  the Voronoi summation formula holds for all $x >0$, namely
$$ \sum_{n=1}^\infty  d(n) f \left(n x\right) -   \int_0^\infty f(xy) (\log y + 2\gamma) dy =  {f(0)\over 4} + 
{1\over x} \sum_{n=1}^\infty  d(n) G \left({n\over x}\right),\eqno(1.40)$$
where
$$ G(x) = \int_0^\infty \left[ 4 K_0\left(4\pi\sqrt{xy}\right) -  2\pi Y_0 \left(4\pi\sqrt{xy}\right)\right] f(y) dy\eqno(1.41)$$
is the integral transform with the combination of the Bessel functions \cite{erd}, Vol. II as the kernel.} 

\begin{proof}  Replacing  $\zeta(s)$ in (1.35) by the right-hand side of  the functional equation (1.4), we  obtain

$$ {1\over 2\pi i}  \int_{\sigma-i\infty}^{\sigma+i\infty}  2^{2s} \pi^{2(s-1)} \zeta^2(1-s) \sin^2\left({\pi s\over 2}\right)\Gamma^2 (1-s) f^*(s) x^{-s} ds$$
$$ = {1\over 2\pi i}  {d\over dx} \int_{1-\sigma-i\infty}^{1-\sigma+i\infty}  2^{2(1-s)} \pi^{-2s} \zeta^2(s) \cos^2\left({\pi s\over 2}\right)\Gamma^2 (s) f^*(1-s) {x^{s}\over s} ds$$ 
$$ = {1\over 2\pi i}  {d\over dx} \int_{1-\sigma-i\infty}^{1-\sigma+i\infty}  2^{1-s} \pi^{-s} \zeta^2(s) \cos\left({\pi s\over 2}\right)\Gamma (s) g^*(1-s) {x^{s}\over s} ds,$$ 
where we denoted by
$$g^*(1-s)= 2^{1-s} \pi^{-s}  \cos\left({\pi s\over 2}\right)\Gamma (s) f^*(1-s)$$
and recalling (1.31),  we find,  correspondingly,
$${1\over x} g\left({1\over x}\right) = \sqrt{2\pi} (F_cf)(2\pi x).\eqno(1.42)$$ 
The differentiation   under the integral sign is possible due to the absolute and uniform convergence and via the asymptotic behavior  (since $x^{-1} (F_cf)(x^{-1}) \in \mathcal{M}_\alpha$)
$$\zeta^2(s) \cos\left({\pi s\over 2}\right)\Gamma (s) g^*(1-s) = O(|t|^{\varepsilon  -3/2}), \  |t| \to \infty,$$ 
where $s= 1-\sigma +it,\ 0 < \sigma < 1.$  The integrand has, possibly, a simple pole at $s=1$ with the residue
$- x f(0)/4$. Moreover, it behaves as $O(|t|^ {-3/2} \log^2 t),\  |t| \to \infty$, when $\sigma=1$ and $O(|t|^{\sigma -5/2}),\ |t| \to \infty$, when $\sigma > 1$.  Hence, moving the contour to the right and  appealing to the residue theorem, the equality (1.35) can be written in the form
$$\sum_{n=1}^\infty  d(n) f \left(n x\right) -   \int_0^\infty f(xy) (\log y + 2\gamma) dy  =  {f(0)\over 4} $$
$$+ {1\over 2\pi i}  {d\over dx} \int_{\sigma-i\infty}^{\sigma+i\infty}  2^{1-s} \pi^{-s} \zeta^2(s) \cos\left({\pi s\over 2}\right)\Gamma (s) g^*(1-s) {x^{s}\over s} ds,\eqno(1.43)$$
where we take again $1 < \sigma <  3/2$ in order to maintain  the absolute convergence of the integral.  Now the latter integral in (1.43) can be treated, employing identity (1.8) and changing the order of integration and summation via the absolute and uniform convergence. Then with the use of (1.31)  we obtain
$${1\over 2\pi i}  {d\over dx} \int_{\sigma-i\infty}^{\sigma+i\infty}  2^{1-s} \pi^{-s} \zeta^2(s) \cos\left({\pi s\over 2}\right)\Gamma (s) g^*(1-s) {x^{s}\over s} ds$$
$$ = {1\over \pi i}  {d\over dx} \sum_{n=1}^\infty  d(n)  \int_{\sigma-i\infty}^{\sigma+i\infty}  \cos\left({\pi s\over 2}\right)\Gamma (s) g^*(1-s) \left({2\pi n \over x}\right)^{-s} {ds\over s}$$
$$=\sqrt{2\pi}   {d\over dx} \sum_{n=1}^\infty  d(n)  \int_0^x   (F_c g)\left({2\pi n\over y}\right) {dy \over y}=  {\sqrt{2\pi} \over x}   \sum_{n=1}^\infty  d(n)  (F_c g)\left({2\pi n\over x}\right),$$
where the term by term differentiation of the series is allowed,  because $g(x)$ satisfies conditions of the theorem due to relation (1.42).    Hence, writing $\sqrt{2\pi}  (F_c g)\left(2\pi n/ x\right)$ in terms of the iterated integral
$$  \sqrt{2\pi} (F_c g)\left({2\pi n\over x}\right) =  4 \int_0^\infty \cos\left({2\pi n  t\over x}\right) {1\over t} 
\int_0^\infty   \cos\left({2\pi  y\over t}\right) f(y) dy dt,$$
and changing formally the order of integration, we invoke the value of the integral (see \cite{prud}, Vol. 1, relation (2.5.24.2)) in terms of the modified Bessel functions
$$\int_0^\infty \cos\left({2\pi n  t\over x}\right)  \cos\left({2\pi  y\over t} \right){dt \over t} = K_0\left(4\pi\sqrt{{ny\over x}}\right) -
{\pi\over 2} Y_0 \left(4\pi\sqrt{{ny\over x}}\right),$$
to arrive at the Voronoi formula (1.40).   To justify this change, we involve integration by parts and some estimates, which are based on the conditions of the theorem.  In fact, we have  for each $x >0$  and $n \in \mathbb{N}$
$$\int_0^\infty \cos\left({2\pi n  t\over x}\right) {1\over t}  \int_0^\infty   \cos\left({2\pi  y\over t}\right) f(y) dy dt
= \left( \int_0^x+ \int_x^\infty\right)  \cos\left({2\pi n  t\over x}\right) {1\over t}  \int_0^\infty   \cos\left({2\pi  y\over t}\right) f(y) dy dt$$
$$=  I_1(x) + I_2(x).$$
Then  by virtue of Fubini's theorem
$$I_1=  - {1\over 2\pi } \int_0^x \cos\left({2\pi n  t\over x}\right)  \int_0^\infty   \sin \left({2\pi  y\over t}\right) f^{(1)} (y) dy dt$$
$$= - {1\over 2\pi } \int_0^\infty  \left(\int_0^x \cos\left({2\pi n  t\over x}\right)   \sin \left({2\pi  y\over t}\right) dt\right)  f^{(1)} (y) dy,$$ 
since 
$$ \int_0^x \left|\cos\left({2\pi n  t\over x}\right) \right| \int_0^\infty  \left| \sin \left({2\pi  y\over t}\right) f^{(1)} (y) dy\right| dt
\le x \int_0^\infty  \left| f^{(1)} (y)\right|  dy  = O(x),$$
and 
$$I_2=    {x\over 2\pi n  } \int_x^\infty  \sin \left({2\pi n  t\over x}\right)  {1\over t^2} \int_0^\infty   \cos \left({2\pi  y\over t}\right) f (y) dy dt$$ 

$$-   {x\over  n  } \int_x^\infty  \sin \left({2\pi n  t\over x}\right)  {1\over t^3} \int_0^\infty   \sin  \left({2\pi  y\over t}\right) f (y) y dy dt $$

$$= {x\over 2\pi n  } \int_0^\infty  \left( \int_x^\infty  \sin \left({2\pi n  t\over x}\right)    \cos \left({2\pi  y\over t}\right) 
{dt \over t^2}\right)   f (y) dy $$ 

$$-   {x\over  n  }  \int_0^\infty  \left( \int_x^\infty  \sin \left({2\pi n  t\over x}\right)   \sin  \left({2\pi  y\over t}\right) 
 {dt \over t^3}\right)   f (y) y dy  $$
because

$${x\over 2\pi n  } \int_0^\infty  \left( \int_x^\infty  \left| \sin \left({2\pi n  t\over x}\right)    \cos \left({2\pi  y\over t}\right) \right| 
{dt \over t^2}\right)  | f (y)| dy \le  {1\over 2\pi n  } \int_0^\infty |f(y)| dy = O(1),$$ 
and

$$  {x\over  n  }  \int_0^\infty  \left( \int_x^\infty  \left| \sin \left({2\pi n  t\over x}\right)   \sin  \left({2\pi  y\over t}\right) \right|
 {dt \over t^3}\right)  | f (y)|  y dy\le   {1\over 2x n  }  \int_0^\infty y |f(y)| dy = O \left(   {1\over x   } \right).$$

\end{proof}

\section{New summation and transformation  formulas}

\subsection{The case $\frac{\zeta(s)}{\zeta(2s)}$.}    The results of this subsection are based on the identities  (1.11), (1.12).  It was shown recently by the author \cite{yakrim}, that they generate for $f\in \mathcal{M}_\alpha$ the so-called reduced M\"{o}bius  operator (cf.  (1.18)).   Namely, we have   
$$(\hat{\Theta} f)(x) =   \sum_{n=1}^\infty |\mu (n)| f(xn) =   \sum_{n: \   \mu(n) \neq 0 }  f(xn)  ,\quad x > 0,\eqno(2.1)$$
where the summation is over all positive integers, which are products of different primes.   Its reciprocal inverse involves  the  Liouville  function  as the kernel 
$$f(x)=  \sum_{n=1}^\infty \lambda(n)   (\hat{\Theta} f)(xn).\eqno(2.2)$$
Moreover, $\hat{\Theta} f$ can be represented by the absolutely convergent integral
$$(\hat{\Theta} f) (x)= {1\over 2\pi i}  \int_{\sigma-i\infty}^{\sigma+i\infty}  \frac{\zeta(s)}{\zeta(2s)}  f^*(s) x^{-s} ds,\  x >0\eqno(2.3)$$
over an arbitrary vertical line in the strip $1 < \sigma < \alpha$.  Hence, as usual, moving the contour to the left and taking into account the analyticity of $f^*(s)$ in the strip  $1/2\le  \sigma < \alpha$,  the residue of the integrand at the simple pole $s=1$ and the absence of zeros of $\zeta(2s)$ on the critical line $s=1/2 +it$ (see \cite{titrim}), we derive the equality 
$$(\hat{\Theta} f) (x)- {6\over \pi^2 x} \int_0^\infty f(y) dy = {1\over 2\pi i}  \int_{\sigma-i\infty}^{\sigma+i\infty}  \frac{\zeta(s)}{\zeta(2s)}  f^*(s) x^{-s} ds,\  x >0,\eqno(2.4)$$
which is valid for $1/2\le  \sigma < 1$. 

{\bf Theorem 4}. {\it Let  $f\in \mathcal{M}_\alpha$. Then the following M\"{u}ntz type formula holds 
$$\frac{\zeta(s)}{\zeta(2s)}  f^*(s) =  \int_0^\infty x^{s-1} \left[ (\hat{\Theta} f) (x)- {6\over \pi^2 x} \int_0^\infty f(y) dy \right] dx,
\eqno(2.5)$$
which is valid for  $1/2\le  \sigma < 1$.}

\begin{proof} Using similar arguments as in the proof of Theorem 2, we recall (1.24) to find 
$$ \frac{\zeta(s)}{\zeta(2s)}  f^*(s) = O( |t|^{\varepsilon - (3+\sigma)/2}), \    1/2\le  \sigma < 1,$$
which guarantees the integrability of the left-hand side in (2.5). Hence the result follows immediately from (2.4) as a reciprocal relation via the Mellin transform (1.1).

\end{proof}

The Mellin transform of $(\hat{\Theta} f)(x) $ exists for $\sigma > 1$ and the process is justifiable as in (1.33) to obtain the equality
$$(\hat{\Theta} f)^*(s) = \frac{\zeta(s)}{\zeta(2s)}  f^*(s),\  1 < \sigma < \alpha.$$
Hence,
$$\frac{(\hat{\Theta} f)^*(s)}{\zeta(s)}  = \frac{f^*(s)}{\zeta(2s)},\    1 < \sigma < \alpha.$$
But the right-hand side of the latter equality is analytic in the strip $1/2 < \sigma < 1$.  Moreover, with the aid of the identity (1.10) and since $f^*(s)$ is integrable, we find the formula of the inverse Mellin transform
$${1\over 2\pi i}  \int_{\sigma-i\infty}^{\sigma+i\infty}  \frac{(\hat{\Theta} f)^*(s)}{\zeta(s)} x^{-s} ds = \sum_{n=1}^\infty  \mu(n)  f(n^2 x),\ x >0,\ 1/2 < \sigma < 1.$$
Besides,   we obtain a sufficient condition for the validity of the Riemann hypothesis.

{\bf Corollary 3.}  {\it Let $f\in  \mathcal{M}_\alpha$ and $(Mf)^*(s)$ is free of zeros in the strip $1/2 < \sigma < 1$.
Then the Riemann hypothesis is true. }

\begin{proof} In fact, since $(\hat{\Theta} f)^*(s) /\zeta(s)$ is analytic in the strip $1/2 < \sigma < 1$ and $(\hat{\Theta} f)^*(s) \neq 0$, it means that $\zeta(s)\neq 0,\  1/2 < \sigma < 1$. Thus possible zeros of the Riemann zeta -function lie only on the line $\sigma = 1/2 $ and the Riemann hypothesis holds true. 
\end{proof}

{\bf Corollary 4}. {\it  Let  $f\in \mathcal{M}_\alpha$. Then the M\"{o}bius  operator $(1.18)$ has the representation in terms of the reduced M\"{o}bius operator $(2.1)$}
$$(\Theta f)(x)=    \sum_{n=1}^\infty   (\hat{\Theta} f) (n^2 x),\quad    x > 0.\eqno(2.6)$$ 

\begin{proof}  We first observe  from (2.1) that $(Mf)(x)$ has the uniform estimate for all $ x >0$
$$ |(\hat{\Theta} f) (x)| \le     \sum_{n=1}^\infty | f(xn) | \le   C x^{-\alpha}   \sum_{n=1}^\infty {1\over n^\alpha}= C_\alpha  \   x^{-\alpha} ,\eqno(2.7)$$
where $C, C_\alpha > 0$ are constants.   Hence, multiplying both sides of (2.5)  by $\zeta(2s)$ and taking the inverse Mellin transform (1.2) for $1/2 < \sigma < 1$, we appeal to (1.21) to derive the equality for the M\"{u}ntz operator
$$(P f)(x)=  {1\over 2\pi i}  \int_{\sigma-i\infty}^{\sigma+i\infty}  \zeta(2s) x^{-s} \int_0^\infty v^{s-1} \left[ (\hat{\Theta} f) (v)- {6\over \pi^2 v} \int_0^\infty f(y) dy \right] dv ds.$$ 
In the meantime, substituting $\zeta(2s)$ in the latter equality by the corresponding series (1.5), one can change the order of  integration and summation owing to (2.5).  Hence we obtain 
$$(P f)(x)=  {1\over 2\pi i}    \sum_{n=1}^\infty \int_{\sigma-i\infty}^{\sigma+i\infty}  (n^2 x) ^{-s} \int_0^\infty v^{s-1} \left[ (\hat{\Theta} f) (v)- {6\over \pi^2 v} \int_0^\infty f(y) dy \right] dv ds,\  1/2 < \sigma < 1.\eqno(2.8)$$ 
But  shifting the contour in the integral of the right-hand side  of (2.4) within the strip $1/2 \le  \sigma < 1$ (this is permitted via analyticity of the integrand in the strip $1/2 < \sigma < 1$ and since it goes to zero when $|t| \to \infty$ uniformly by $\sigma$ in each inner strip), we verify the integrability of its left-hand side over $\mathbb{R}_+$ with respect to the measure $x^{\sigma -1} dx,\ 1/2 < \sigma < 1$, i.e.
$$h(x)= (\hat{\Theta} f) (x)- {6\over \pi^2 x} \int_0^\infty f(y) dy \in L_1( \mathbb{R}_+; x^{\sigma -1} dx).$$
Further, its Mellin transform (1.1) $h^*(s)$  is integrable due to (2.5). Hence (cf. \cite{tit}) formula (2.8) can be simplified and we get
$$(P f)(x)=    \sum_{n=1}^\infty  \left[ (\hat{\Theta} f) (n^2 x)- {6\over (\pi  n)^2 x} \int_0^\infty f(y) dy \right],\  x >0.$$ 
Now, splitting in two series, which is possible by virtue of the estimate (2.7), and taking into account (1.18), (1.20) and the value of $\zeta(2)$, we established  (2.6)  and completed the proof of Corollary 4. 
\end{proof}

Analogously, defining the operator with the Liouville function 
$$(\Lambda f)(x)=    \sum_{n=1}^\infty   \lambda (n) f(n x),\quad    x > 0,\eqno(2.9)$$ 
and appealing to (1.12), we find
$$(\Lambda f)(x)= {1\over 2\pi i}  \int_{\sigma-i\infty}^{\sigma+i\infty}  \frac{\zeta(2s)}{\zeta(s)}  f^*(s) x^{-s} ds,\  x >0\eqno(2.10)$$
when $\sigma > 1$. Reciprocally,
$$\zeta(s) (\Lambda f)^*(s) =  \zeta(2s)  f^*(s),\  \sigma > 1\eqno(2.11)$$
and  therefore
$$ \sum_{n=1}^\infty   f(n^2 x) =  {1\over 2\pi i}  \int_{\sigma-i\infty}^{\sigma+i\infty}  \zeta(s)   (\Lambda f)^*(s)  x^{-s} ds,\  x >0, \ \sigma > 1.$$
Hence, replacing $\zeta(s)$ in the latter integral by its series (1.5), we derive
$$  \sum_{n=1}^\infty   f(n^2 x) =  {1\over 2\pi i}   \sum_{n=1}^\infty \int_{\sigma-i\infty}^{\sigma+i\infty}   (\Lambda f)^*(s)  (nx)^{-s} ds,\  x >0, \   \sigma > 1,\eqno(2.12)$$
where the change of the order of integration and summation is allowed via the absolute convergence (see (2.11)). 
Similarly as above from (2.9) we have
$$|(\Lambda f)(x) | \le  \sum_{n=1}^\infty |\lambda(n)  f(n x) | \le \sum_{n=1}^\infty |  f(n x) | = O(x^{-\alpha}),\  x >0$$
and  one shows from (2.10) that $(\Lambda f)(x) \in L_1( \mathbb{R}_+; x^{\sigma -1} dx),\   1 < \sigma <  \alpha$.  Hence (2.12)  becomes 
$$  \sum_{n=1}^\infty   f(n^2 x) =    \sum_{n=1}^\infty  (\Lambda f)(nx),\  x >0\eqno(2.13)$$
 and we proved the following 

 {\bf Theorem 5}.  {\it Let $f \in \mathcal{M}_\alpha$. Then for all $x >0$ the summation formula $(2.13)$ takes place for the operator  with the Liouville function $(2.9)$.}
 
 A necessary condition for the validity of the Riemann hypothesis is given by 
 
 {\bf Corollary 5}. {\it Let the Riemann hypothesis holds true. Then for any  $f \in \mathcal{M}_\alpha$ the summation formula takes place for all $ x>0$}
 $$ \sum_{n=1}^\infty   \lambda (n) f(n x)  =   \sum_{n=1}^\infty  {1\over 2\pi i}  \int_{\sigma-i\infty}^{\sigma+i\infty}  \frac{f^*(s)}{\zeta(s)}   (n^2x)^{-s} ds,\  1/2 < \sigma < 1.$$
 
 \begin{proof}  The proof follows immediately from (2.9), (2.10), the analyticity of the integrand in (2.10) in the strip 
$1/2 < \sigma < 1$ under the truth of the Riemann hypothesis and its integrability via the asymptotic behavior
(see \cite{titrim}, Ch. XIV)
$$\frac{f^*(s)}{\zeta(s)} = O( |t|^{\varepsilon - 2}),\ |t| \to \infty,\   1/2 < \sigma < 1.$$ 
 
 \end{proof}

 Meanwhile, formula (2.4) will be a starting point to prove  summation formulas  for the operator (2.1). Precisely, it states by

  {\bf Theorem 6}. {\it Let $f \in \mathcal{M}_\alpha$. Then for all $x >0$ the following Poisson type summation  formulas  hold
  $$  \sum_{n:\  \mu(n)\neq 0} f(xn)  =   \sum_{n,m =1}^\infty \mu(n)   f (n^2 m x)
  = \sum_{n=1}^\infty \mu(n)  (\Theta f) (n^2 x),\eqno(2.14)$$
 $$  \sum_{n:\  \mu(n)\neq 0} f(xn)   =  {3 \sqrt 2\over  x \pi\sqrt \pi} (F_cf)(0)  
+  \sum_{n =1}^\infty \mu(n) \left[   {\sqrt{2\pi} \over n^2 x}   \sum_{m =1}^\infty (F_cf)\left({2\pi m\over n^2 x}\right)-   {1\over 2} f (0)\right].\eqno(2.15)$$
Finally, if, in addition, $f(x)= O(x^\beta),\  \beta > 1/2, \ x \to 0$,  then  the summation formula takes place 
 
 $$  \sum_{n:\  \mu(n)\neq 0}  f(xn) - {6\over \pi^2 x} \int_0^\infty f(y) dy = 2^{-3/2} e^{-i\pi/4} \sqrt x  \sum_{n,m  =1}^\infty    {\mu(n)\over n^2 m^{3/2} } \   G\left( {\pi x \over 2n^2m}\right),\eqno(2.16)$$
where
$$G(x)=  \int_0^\infty { f(1/ u) \over \sqrt u }  \left[ e^{i x u} \hbox{erf} \left( e^{i\pi/4}\sqrt{ x u } \right) + e^{- i x u} \hbox{erfi} \left( e^{i\pi/4}\sqrt{ x u }\right)\right]   du\eqno(2.17) $$
is the Mellin convolution type transform with  a combination of the error functions \cite{prud}, Vol. 2  as the kernel.}

\begin{proof}   In fact, in order to prove equalities (2.14), we appeal to the equalities (1.18), (1.20), (1.21) and the Ramanujan identity (1.10).  Hence,  substituting the value $[\zeta(2s)]^{-1},\  \sigma > 1/2$ into the right- hand side of (2.4), we  change  the order of integration and summation by the absolute convergence.   Thus
$${1\over 2\pi i}  \int_{\sigma-i\infty}^{\sigma+i\infty}  \frac{\zeta(s)}{\zeta(2s)}  f^*(s) x^{-s} ds = 
 \sum_{n =1}^\infty \mu(n) (Pf)( n^2 x) $$
 $$=  \sum_{n =1}^\infty \mu(n) \left[   (\Theta f) (n^2 x) - {1\over n^2 x} \int_0^\infty f(y) dy \right] 
 = \sum_{n =1}^\infty \mu(n) (\Theta f) (n^2 x) - {6\over \pi^2 x} \int_0^\infty f(y) dy, $$
where the splitting in two series is possible due to the definition of the class  $\mathcal{M}_\alpha$ and their absolute convergence.  Hence combining with (2.4), we come up with equalities (2.14).   Similar arguments are used to prove (2.15), where the application of the Poisson formula (1.27) in the class $\mathcal{M}_\alpha$ is involved (see the proof of Theorem 1).   In fact,  for each $n \in \mathbb{N}$ and $ x > 0$,  we have 
$$ \sum_{m=1}^\infty   f (n^2m x)=  {\sqrt{\pi} \over n^2 x\sqrt 2 }  (F_cf) (0) +  {\sqrt{2\pi} \over n^2 x} \sum_{m =1}^\infty   (F_cf) \left({2\pi m\over n^2 x}\right)-  {1\over 2} f (0). $$
Therefore,
$$ \sum_{n,m =1}^\infty \mu(n)   f (n^2 m x) =   {3 \sqrt 2\over  x \pi\sqrt \pi} (F_cf)(0)  
+  \sum_{n =1}^\infty \mu(n) \left[   {\sqrt{2\pi} \over n^2 x}   \sum_{m =1}^\infty (F_cf)\left({2\pi m\over n^2 x}\right)-   {1\over 2} f (0)\right],$$
where the splitting of the series is possible via them convergence and identity (1.10) is applied.   Combining with (2.14),  it  gives (2.15).

Let us prove (2.16).  Applying  the functional equation (1.4), replacing  $s$ by $1-s$ and performing the differentiation under the integral sign, we move the contour to the right in the obtained integral,  since $f(0)=0$ and therefore the integrand has a removable singularity at the point $s=1$. Hence we   write the right-hand side of  (2.4)  as follows 
$$ {1\over 2\pi i}  {d\over dx} \int_{\sigma-i\infty}^{ \sigma+i\infty}  \frac{\zeta(s)}{\zeta(2(1-s))}   2^{1-s} \pi^{-s}   \cos\left({\pi s\over 2}\right)\Gamma(s) f^*(1-s) x^{s} {ds\over s},\ 1 < \sigma < \hbox{min }(\alpha,\ 3/2). $$
Hence applying again  the functional equation (1.4) in the denominator  and employing the supplement  and duplication formulas  for the gamma-function \cite{erd}, Vol. I,  we obtain
$${1\over 2\pi i}  {d\over dx} \int_{\sigma-i\infty}^{ \sigma+i\infty}  \frac{\zeta(s)}{\zeta(2(1-s))}   2^{1-s} \pi^{-s}   \cos\left({\pi s\over 2}\right)\Gamma(s) f^*(1-s) x^{s} {ds\over s}$$
$$ =  {1\over 2\pi i}  {d\over dx}   \int_{\sigma-i\infty}^{ \sigma+i\infty}  \frac{\zeta(s)}{\zeta(2s-1)}   \frac{2^{3/2-2s} \pi^{s- 1}\Gamma(s/2)\Gamma(1-s/2)  }{  \Gamma((s-1/2)/2) \Gamma((s +1/2)/2) } f^*(1-s) x^{s} {ds\over s},  \  1 < \sigma < \hbox{min }(\alpha,\ 3/2).$$
Further, the series representations (1.5) and (1.10) and the possibility to change  the order of integration and summations drive us at the equality
$${1\over 2\pi i}  {d\over dx}   \int_{\sigma-i\infty}^{ \sigma+i\infty}  \frac{\zeta(s)}{\zeta(2s-1)}   \frac{2^{3/2-2s} \pi^{s- 1}\Gamma(s/2)\Gamma(1-s/2)  }{  \Gamma((s-1/2)/2) \Gamma((s +1/2)/2) } f^*(1-s) x^{s}{ ds\over s}$$
$$=  {2^{3/2} \over \pi }   {d\over dx}  \sum_{n,m  =1}^\infty    { n \mu(n) \over 2\pi i}   \int_{\sigma-i\infty}^{ \sigma+i\infty}  \frac {\Gamma(s/2)\Gamma(1-s/2)  }{  \Gamma((s-1/2)/2) \Gamma((s +1/2)/2) } f^*(1-s) \left({4n^2 m\over \pi x}\right)^{-s} {ds\over s}.$$
But since $f(0)=0$, one can move the contour to the left, considering the latter integral for $1/2 < \sigma < 1$.   Hence we take  into account the asymptotic behavior
$$   \frac{\Gamma(s/2)\Gamma(1-s/2)  }{  s\   \Gamma((s-1/2)/2) \Gamma((s +1/2)/2) }  
= O(|t|^{-\sigma}),\ |t| \to \infty, $$
to  recall the Parseval equality (1.3) and the value of the integral
$${1 \over 2\pi i} \int_{\sigma-i\infty}^{ \sigma+i\infty}  \frac {\Gamma(s/2)\Gamma(1-s/2)  }{  \Gamma((s-1/2)/2) \Gamma((s +1/2)/2) } \   u^{-s}\   {ds\over s}$$
$$=  {e^{-i\pi/4} \over \sqrt \pi} \int_0^{1/u}  \sqrt  y \left[ e^{2iy} \hbox{erf} \left( e^{i\pi/4}\sqrt{2y}\right)+ e^{-2iy} \hbox{erfi} \left( e^{i\pi/4}\sqrt{2y}\right)\right] dy,$$
which is calculated, in turn,  by the Slater theorem and relation (7.14.2.75) in  \cite{prud}, Vol. 3 and  contains  a combination of the error functions.    Therefore we  get the equality 
$${2^{3/2} \over \pi }   {d\over dx}  \sum_{n,m  =1}^\infty    { n \mu(n) \over 2\pi i}   \int_{\sigma-i\infty}^{ \sigma+i\infty}  \frac {\Gamma(s/2)\Gamma(1-s/2)  }{  \Gamma((s-1/2)/2) \Gamma((s +1/2)/2) } f^*(1-s) \left({4n^2 m\over \pi x}\right)^{-s} {ds\over s}$$
$$ =  {2^{3/2} e^{-i\pi/4} \over \pi\sqrt \pi  }   {d\over dx}  \sum_{n,m  =1}^\infty    n \mu(n) \int_0^\infty f(u) \int_0^{\pi x/(4 n^2 m\  u)}   \sqrt  y \left[ e^{2iy} \hbox{erf} \left( e^{i\pi/4}\sqrt{2y}\right)\right.$$
$$\left. + e^{-2iy} \hbox{erfi} \left( e^{i\pi/4}\sqrt{2y}\right)\right]  dy\  du. \eqno(2.18)$$ 
Formally, performing the differentiation term by term in the double series and under the integral sign  in the right-hand side of the latter equality,  we find it in the form 
$$ 2^{-3/2} e^{-i\pi/4} \sqrt x  \sum_{n,m  =1}^\infty    {\mu(n)\over n^2 m^{3/2} } \int_0^\infty { f(1/ u) \over \sqrt u }  
\left[ e^{i\pi x u /(2n^2 m )} \hbox{erf} \left( e^{i\pi/4}\sqrt{{\pi x u \over 2n^2 m } }\right)\right.$$
$$\left. + e^{- i\pi x u /(2n^2 m )} \hbox{erfi} \left( e^{i\pi/4}\sqrt{{\pi x u \over 2n^2 m }}\right)\right]   du.\eqno(2.19) $$ 
In order to justify this operation, we appeal to the definition of the error functions \cite{prud}, Vol.  2 and write the expression in the square brackets as 
$$ 2^{-3/2} e^{-i\pi/4} \sqrt x  \left[ e^{i\pi x u /(2n^2 m )} \hbox{erf} \left( e^{i\pi/4}\sqrt{{\pi x u \over 2n^2 m } }\right)+ e^{- i\pi x u /(2n^2 m )} \hbox{erfi} \left( e^{i\pi/4}\sqrt{{\pi x u \over 2n^2 m }}\right)\right] $$
$$ =   {x\over n} \sqrt{{ u \over  m }} \int_0^1 \cos\left( {\pi x u (1-t^2)  \over 2n^2 m }\right) dt.$$
But for $x >0$
$$\left| \int_0^1 \cos\left(x (1-t^2) \right) dt\right| = \left|\cos x  \int_0^1 \cos\left(x t^2 \right) dt + \sin x  \int_0^1 \sin\left(x t^2 \right) dt\right| $$
$$=   \left|{\cos x \over \sqrt x} \left[  \int_0^\infty \cos\left( t^2 \right) dt - \int_{\sqrt x}^\infty \cos(t^2) dt \right]    +  {\sin x\over \sqrt x }  \left[ \int_0^\infty \sin \left( t^2 \right) dt - \int_{\sqrt x}^\infty \sin (t^2) dt \right]  \right| $$
$$ \le {1\over \sqrt x} \left[ \left| \int_0^\infty \cos\left( t^2 \right) dt \right| + \left| \int_{\sqrt x}^\infty \cos(t^2) dt \right|     +   \left| \int_0^\infty \sin \left( t^2 \right) dt\right| + \left| \int_{\sqrt x}^\infty \sin (t^2) dt \right|  \right] \le {C\over \sqrt x},$$ 
where $C > 0$ is an absolute constant.   Hence under conditions on the theorem 
$$ \sum_{n,m  =1}^\infty    {|\mu(n)|\over n^2 m^{3/2} } \int_0^\infty { |f(1/ u)|  \over \sqrt u }  
\left |e^{i\pi x u /(2n^2 m )} \hbox{erf} \left( e^{i\pi/4}\sqrt{{\pi x u \over 2n^2 m } }\right)+ e^{- i\pi x u /(2n^2 m )} \hbox{erfi} \left( e^{i\pi/4}\sqrt{{\pi x u \over 2n^2 m }}\right)\right|   du $$
$$\le C_1  \sum_{n  =1}^\infty    {1\over n^2 }  \sum_{m  =1}^\infty    {1\over m^{3/2} } \left[ \int_0^1 O(u^{\alpha -1/2}) du +
  \int_1^\infty  O(u^{-\beta -1/2}) du \right]  < \infty,  $$ 
where $C_1 >0$ is a constant,  and the differentiation is allowed in (2.18).   Thus returning to (2.19) and combining with (2.4), we arrive at the equality (2.16), completing the proof of Theorem 6. 

\end{proof}

\subsection{The case ${\zeta(s-1)\over \zeta(s)}$.}   Employing identities  (1.15),  (1.16) we introduce the transformations, involving Euler's  totient function $\varphi(n)$ and divisor function $a(n)$, respectively,
$$(\Phi f)(x)= \sum_{n=1}^\infty  \varphi(n) f(xn),\quad x > 0,\eqno(2.20)$$
$$(A f)(x)= \sum_{n=1}^\infty  a(n) f(xn),\quad x > 0.\eqno(2.21)$$
The M\"{u}ntz type formulas for these operators  can be established in the same manner as in Theorem 4 and we have

{\bf Theorem 7}. {\it Let  $f\in \mathcal{M}_\alpha$. Then the following M\"{u}ntz type formulas  hold 
$$\frac{\zeta(s-1)}{\zeta(s)}  f^*(s) =  \int_0^\infty x^{s-1} \left[ (\Phi f)(x)- {6\over \pi^2 x} \int_0^\infty f(y) dy \right] dx,
\eqno(2.22)$$
$$\frac{1- 2^{1-s}}{1-2^{-s}}  \zeta(s-1) f^*(s) =  \int_0^\infty x^{s-1} \left[ (A f)(x)- {2\over 3 x} \int_0^\infty f(y) dy \right] dx,
\eqno(2.23)$$
which are  valid for  $1 \le  \sigma < 2$.}

Now, for $\sigma > 2$, we have (cf. \cite{yakrim})
$$\zeta(s-1) f^*(s) = \zeta(s) (\Phi f)^*(s),$$
$$(1- 2^{1-s})  \zeta(s-1) f^*(s) =  (1-2^{-s}) (A f)^*(s).$$
Hence,  taking into account the asymptotic behavior of the totient function (see above) and identity (1.16) with the divisor function $a(n)$, we get

 {\bf Theorem 8}.  {\it Let $f \in \mathcal{M}_\alpha$ with $\alpha > 2$. Then for all $x >0$ the following summation formulas  hold valid }
$$  \sum_{n=1}^\infty  n  f(xn) =    \sum_{n=1}^\infty  (\Phi f)(nx),\eqno(2.24)$$
$$  (A f)(x)-   (A f)(2x) =  \sum_{n=1}^\infty  n \left[ f(xn)  - 2 f(2xn) \right] =  \sum_{n=1}^\infty  \left[ (\Phi f)(nx)-  2 (\Phi f)(2 nx)\right].\eqno(2.25)$$

Further, the M\"{u}ntz formula (1.19),  Lemma 1 and Theorem 7  lead us to 

{\bf Theorem 9}.  {\it Let $f \in \mathcal{M}_\alpha$ with $\alpha > 2$. Then for all $x >0$ the  summation formulas 
$$\sum_{n=1}^\infty  \varphi(n) f(xn) -   \sum_{n,m  =1}^\infty    m \mu(n) f(n m x) = 
{6\over \pi^2} \int_0^\infty (1-y) f(xy) dy,\eqno(2.26)$$
$$\sum_{n=1}^\infty  \left[ a(n)- n\right]  f(xn)    +   \sum_{n,m =1}^\infty    n f( 2^m x  n) = 
 {2\over 3} \int_0^\infty (1-y) f(xy) dy \eqno(2.27)$$
 take place. }
 
 {\begin{proof}  In fact, as an immediate consequence of the M\"{u}ntz type formula (2.22), identity (1.10) the asymptotic behavior 
 $$\frac{\zeta(s-1)}{\zeta(s)}  f^*(s) = O(|t|^{\varepsilon -1- \sigma/2}),\   |t| \to \infty,\ 1 < \sigma < 2$$
 and the inversion formula (1.2) for the Mellin transform, we have  the chain of equalities 
 $$ \sum_{n=1}^\infty  \varphi(n) f(xn) -  {6\over \pi^2 x} \int_0^\infty f(y) dy = {1\over 2\pi i} \sum_{n=1}^\infty  \mu(n) 
  \int_{\sigma-1 -i\infty}^{ \sigma-1 +i\infty}   \zeta(s) f^*(1+s) (xn)^{-s -1} ds$$ 
 $$=  \sum_{n =1}^\infty  \mu(n)  \left[   \sum_{m =1}^\infty  m f(nm x) -  {1\over (xn)^2} \int_0^\infty y f(y) dy\right]
 = \sum_{n, m  =1}^\infty  m \mu(n) \  f(nm x) -  {6\over (\pi x )^2} \int_0^\infty y f(y) dy. $$
 Hence we easily come up with (2.26).  Similarly,
 $$\sum_{n=1}^\infty  a(n) f(xn)  - {2\over 3 x} \int_0^\infty f(y) dy =  {1\over 2\pi i}  
  \int_{\sigma-1 -i\infty}^{ \sigma-1 +i\infty}   \zeta(s) f^*(1+s) x^{-s -1} ds$$
  $$ -  \sum_{m =1}^\infty   \int_{\sigma-1 -i\infty}^{ \sigma-1 +i\infty}   \zeta(s) f^*(1+s) \left(2^m\ x \right) ^{-s -1} ds 
  =   \sum_{m =1}^\infty  m f(m x) - {1\over x^2} \int_0^\infty y f(y) dy  $$
 $$-  \sum_{m =1}^\infty  \left[ \sum_{n =1}^\infty  n f( 2^m\ x \ n) -   {2^{-2m} \over x^2} \int_0^\infty y f(y) dy \right] 
 = \sum_{m =1}^\infty  m f(m x) - \sum_{n,m =1}^\infty    n f( 2^m\ x \ n) $$ 
 $$-  {2\over 3 x^2} \int_0^\infty y f(y) dy,$$  
 which gives (2.27). 
 
 \end{proof} 
 
 \subsection{The generalized Voronoi operator.}    Let us consider the $k$-th iteration of the M\"{u}ntz operator (1.20)
 $(P^{k}f)(x),\ k \in \mathbb{N}_0,\  (P^{0}f)(x) \equiv f(x)$, assuming the  conditions $(P^{j}f)(x) \in  \mathcal{M}_\alpha
 ,\ j=0,1,\dots k-1,\  k \ge 1.$  Then similar to (1.37) we define the generalized Voronoi operator $(\mathbb{V}_k f)(x)$ as
 $$  (\mathbb{V}_k f)(x)=  (P^{k}f)(x),\ x >0,\eqno(2.28)$$
 letting   $(\mathbb{V}_2 f)(x) \equiv  (Vf)(x)$ via Corrolary 2.   Further, taking identity (1.9), we derive analogously to (1.32), (1.33), (1.34) the following relations
 $$  (D_k f)(x) =  \sum_{n=1}^\infty  d_k(n) f \left(n x\right) = {1\over 2\pi i}  \int_{\sigma-i\infty}^{\sigma+i\infty}  \zeta^k(s) f^*(s) x^{-s} ds,\ x >0, \eqno(2.29)$$ 
$$(D_kf)^*(s)= \zeta^k(s) f^*(s),\eqno(2.30)$$
where $ 1 < \sigma < \alpha$.  On the other hand,  $ \zeta^k(s) f^*(s)$ is analytic in the strip $0 < \sigma < \alpha$ except $s=1$, where there is a pole of order $k$.  Moreover, as above
$$\zeta^k(s) f^*(s) = O( |t|^{k(\varepsilon+ (1-\sigma)/2)-2 }),\ |t| \to \infty,\ 0 < \sigma < 1,$$
$$\zeta^k(s) f^*(s) = O( |t|^{- 2}\log^k t),\ |t| \to \infty,\   \sigma = 1.$$
Therefore, when $\hbox{max} \left(0, 1- {2\over k}\right)  < \sigma < 1$, it guarantees the integrability in this strip and possibility to move the contour to the left in the integral (2.29), counting the residue at the multiple pole $s=1$.  Then since 
$${\rm Res}_{s=1} [ \zeta^k(s) f^*(s) x^{-s}] =   {1\over (k-1) !} \lim_{s\to 1} {d^{k-1}\over s^{k-1}} \left[ 
\left( (s-1) \zeta(s)\right)^k  f^*(s) x^{-s} \right]$$
$$=  {1\over (k-1) !} \lim_{s\to 1} \sum_{r=0}^{k-1}  { k-1 \choose r } \  \left[\left( (s-1) \zeta(s)\right)^k  \right]^{(r)} \left[ f^*(s) x^{-s} \right] ^{(k-1-r)}$$
and, in turn, 
$$\left[ f^*(s) x^{-s} \right] ^{(k-1-r)} = \sum_{m =0}^{k-1-r}  { k-1-r \choose m }   \left( f^*(s)\right)^{(k-1-r-m)} \left( x^{-s} \right)^{(m)} $$
$$ =  x^{-s} \int_0^\infty f(y) y^{s-1}   \sum_{m =0}^{k-1-r}  { k-1-r \choose m }  (-1)^m (\log x)^m  (\log y)^{k-1-r-m} dy$$
$$=  x^{-s} \int_0^\infty f(y) y^{s-1}  \left( \log \left({y\over x} \right) \right) ^{k-1-r}  dy =  \int_0^\infty f(x y) y^{s-1} 
 \left( \log y \right) ^{k-1-r}  dy.$$
Meanwhile,  in order to calculate the $r$-th derivative of $\left((s-1) \zeta(s)\right)^k$ we appeal to the familiar  Fa\' {a} di Bruno formula \cite{bruno}.   Thus we obtain 
$$ \left[ \left((s-1)\zeta(s)\right)^k \right] ^{(r)} =  \sum { r! \  k! \  \left( (s-1) \zeta(s)\right)^{k-n}  \over  (k-n)!  b_1! b_2!\dots b_r! }  \left( {\left((s-1) \zeta(s)\right)^{(1)} \over 1!} \right) ^{b_1} \dots  \left( {\left((s-1) \zeta(s)\right) ^{(r)}) \over r! }\right) ^{b_r},$$ 
where the sum is over all different solutions in nonnegative integers $b_1,b_2,\dots, b_r$ of  
$b_1+2b_2+\dots+ rb_r =r$ and $ n= b_1+b_2+\dots + b_r$.   On the other hand, since the Laurent series of the Riemann zeta-function in the neighborhood of $s=1$ has the form (see \cite{erd}, Vol. I)
$$\zeta(s)= {1\over s-1} + \gamma + \sum_{m=1}^\infty \gamma_m(s-1)^m,$$
where $\gamma $ is the Euler constant and
$$\gamma_m= {(-1)^m\over m!} \lim_{l\to \infty} \left[\sum_{j=1}^l j^{-1} \log^m j - (m+1)^{-1} \log^{m+1} l \right]$$
are Stieltjes constants, we easily observe that 
$$\lim_{s\to 1} \left((s-1) \zeta(s)\right)^{(1)} = \gamma,$$
 $$\lim_{s\to 1} \left((s-1) \zeta(s)\right)^{(m)} = m! \gamma_{m-1},\ m = 2,3,\dots, r.$$
 Thus combining with the above calculations we finally obtain
 $${\rm Res}_{s=1} [ \zeta^k(s) f^*(s) x^{-s}] =    \int_0^\infty f(x y) P_{k-1} (\log y) dy,$$
 where $P_{k-1}(x)$ is a polynomial of degree $k-1$ and we give its explicit form, which seems to be new (cf. \cite{titrim}, p. 313), namely
 $$P_{k-1} (x)=  x^{k-1}  +  \sum_{r=1}^{k-1} c_{k,r} x^{k-1-r},\eqno(2.31)$$
 $$c_{k,r} =   {k! \over ( k-1-r)!} \ \sum { \gamma^{b_1}  \gamma^{b_2}_1 \dots  \gamma_{r-1}^{b_r}  \over  (k-n)!  b_1! b_2!\dots b_r! },\eqno(2.32)$$
 where  $\ \gamma_0\equiv \gamma$ and the latter sum is,  as above,  over all different solutions in nonnegative integers $b_1,b_2,\dots, b_r$ of  $b_1+2b_2+\dots+ rb_r =r$ and $ n= b_1+b_2+\dots + b_r$. In particular, letting $k=2$, we get immediately the  residue in the case of Voronoi's  operator (see (1.35)).  So, returning to (2.28), (2.29),  we come up with the representation of the generalized Voronoi operator 
$$(\mathbb{V}_k f)(x)= {1\over 2\pi i}  \int_{\sigma-i\infty}^{\sigma+i\infty}  \zeta^k(s) f^*(s) x^{-s} ds =  (D_k f)(x)
-   \int_0^\infty f(x y) P_{k-1} (\log y) dy ,\  x >0,\eqno(2.33)$$
where $ \hbox{max} \left(0, 1- {2\over k}\right)  < \sigma < 1,\ k \ge 1$.   An analog of Theorem 2 is

{\bf Theorem 10}. {\it Let $f \in \mathcal{M}_\alpha,\ k \in \mathbb{N}$. Then the M\"{u}ntz type formula for the generalized Voronoi operator  $(2.33)$
$$\zeta^k(s) f^*(s) = \int_0^\infty  x^{s-1} \left[  \sum_{n=1}^\infty  d_k(n) f \left(n x\right) -   \int_0^\infty f(xy)  P_{k-1} (\log y)dy\right] dx\eqno(2.34)$$
is valid for $\hbox{max} \left(0, 1- {2\over k}\right)   < \sigma < 1$, where the polynomial $P_{k-1}(x)$ is defined by $(2.31), (2.32)$.} 

{\bf Remark 1}.  For the classical M\"{u}ntz operator (1.20) ($k=1$) the corresponding polynomial is $P_0(x) \equiv 1$.

\subsection{The case ${\zeta^{k+1} (s)\over \zeta(2s)},\   k\in \mathbb{N}$.}   This case is devoted to formulas, involving the composition of the generalized Voronoi operator (2.28) and reduced M\"{o}bius operator (2.1).  Precisely, basing on the M\"{u}ntz  formula  (1.19),  Voronoi formula (1.35)  and  M\"{u}ntz  type formula  ( (2.34),  one can prove in  the same manner the following

{\bf Theorem 11}. {\it Let $f, \  (\hat{\Theta} f) (x)  \in \mathcal{M}_\alpha,\ k \in \mathbb{N}$. Then the M\"{u}ntz type formula 
$${\zeta^{k+1} (s)\over \zeta(2s)}  f^*(s) = \int_0^\infty  x^{s-1} \left[  \sum_{n=1}^\infty  d_k(n)  (\hat{\Theta} f)(n x) -   \int_0^\infty (\hat{\Theta} f)(xy)   P_{k-1} (\log y)dy\right] dx$$
is valid for $\hbox{max} \left(1/2, 1- {2\over k}\right)   \le  \sigma < 1$.  Reciprocally,
$$ \sum_{n=1}^\infty  d_k(n)  (\hat{\Theta} f)(n x) -   \int_0^\infty (\hat{\Theta} f)(xy)   P_{k-1} (\log y)dy
=  {1\over 2\pi i}  \int_{\sigma-i\infty}^{\sigma+i\infty}   {\zeta^{k+1} (s)\over \zeta(2s)} f^*(s) x^{-s} ds.$$
In particular,  employing identities $(1.7), (1.13), (1.14)$, the following summation formulas with arithmetic functions $\omega(n)$ and $d(n)$ hold for all $x >0$}
$$ \sum_{n=1}^\infty  2^{\omega(n)} f (n x)  = \sum_{n=1}^\infty   (\hat{\Theta} f)(n x) -   \int_0^\infty (\hat{\Theta} f)(xy) dy,$$
$$ \sum_{n=1}^\infty  d(n^2)  f (n x) =  \sum_{n=1}^\infty   d(n) (\hat{\Theta} f)(n x) -   \int_0^\infty (\hat{\Theta} f)(xy) (\log y  + 2\gamma) dy,$$
$$ \sum_{n=1}^\infty  d^2(n)  f (n x) =  \sum_{n=1}^\infty   d_3(n) (\hat{\Theta} f)(n x) -   \int_0^\infty (\hat{\Theta} f)(xy) \left[\log^2 y + 3\gamma(\log y +2\gamma) + 3\gamma_1\right] dy.$$

{\bf Remark 2}. Analogously,  one can deduce summation formulas, which are associated with identities (1.7), (1.13), (1.14).  It will contain  integral transforms of the Mellin convolution type with the hypergeometric functions ${}_3F_2,\   {}_5F_2$ and ${}_7F_2$ \  (cf. \cite{erd}, Vol. I,\  \cite{yal}) , respectively.  We leave this task to the interested reader.

\section{Particular examples}

In this section we will consider curious particular examples of the above summation formulas, involving series with the summation over positive integers, which are products of different primes or over positive integers, containing any factor of power higher than the first.  Let, for instance, $f(x)= e^{-x} \in    \mathcal{M}_\alpha$.  Then after simple calculation of its Fourier cosine transform,   formulas (2.14), (2.15)  yield the following identities
$$  \sum_{n:\  \mu(n)\neq 0} e^{-nx}  =   \sum_{n =1}^\infty {\mu(n)\over e^{n^2 x}-1} 
=    {6\over  x \pi^2} +  \sum_{n =1}^\infty \mu(n) \left[   \sum_{m =1}^\infty {2n^2 x\over 4\pi^2 m^2+ n^4 x^2} 
-   {1\over 2} \right],\ x >0.\eqno(3.1)$$ 
Meanwhile,    (see relation (5.4.5.1) in \cite{prud} )
$$ \sum_{m =1}^\infty {2x\over m^2+ x^2} =   \pi \coth\left(\pi x\right)  - {1\over x}.$$
So, the second equality in  (3.1) becomes 
$$   \sum_{n =1}^\infty {\mu(n)\over e^{n^2 x}-1} 
=    {6\over  x \pi^2} +  {1\over 2} \sum_{n =1}^\infty {\mu(n) \over \cosh(n^2x)},\ x >0\eqno(3.2)$$ 
or,
$$ \sum_{n =1}^\infty \mu(n)\left[ {1\over e^{n^2 x}-1}-   {1\over 2\cosh(n^2x)} \right]=\    {6\over  x \pi^2} ,\ x >0.$$
On the other hand, taking the Fourier sine transform of both sides of the first equality in (3.1), we appeal to relations
(2.5.30.8), (2.5.34.4) in \cite{prud}, Vol. 1 to derive
$$ \sum_{n:\  \mu(n)\neq 0} {2x\over n^2+ x^2}  =   \sum_{n =1}^\infty \mu(n)\left[ {\pi\over n^2 } \coth\left({\pi x \over n^2}\right)  - {1\over x} \right] ,\ x >0.\eqno(3.3)$$
Further, integrating with respect to $x$ in (3.3), it gives another curious formula
$$ \sum_{n:\  \mu(n)\neq 0} \log\left({n^2+ x^2\over n^2+1} \right)   =   \sum_{n =1}^\infty \mu(n)
\log\left({ \sinh (\pi x/n^2)\over x \sinh (\pi/n^2)} \right) ,\ x >0.$$ 
In particular, when $x \to 0$, we find
$$ \sum_{n:\  \mu(n)\neq 0}\log\left({n^2\over n^2+1} \right)   =   \sum_{n =1}^\infty \mu(n)
\log\left({ \pi \over n^2 \sinh (\pi/n^2)} \right).$$ 
In the meantime,   
$$\sum_{n =1}^\infty e^{-nx} = {1\over e^x-1}.$$
Hence we get immediately from (3.1), (3.2)
$$\sum_{n:\  \mu(n)= 0} e^{-nx}  =  \sum_{n =2}^\infty {\mu(n)\over 1- e^{n^2 x}}
= {1\over e^x-1} -  {6\over  x \pi^2} - {1\over 2} \sum_{n =1}^\infty {\mu(n) \over \cosh(n^2x)}.\eqno(3.4)$$ 

Finally,   we will consider the case of the self-reciprocal Fourier cosine transforms,  and certainly the most important and  familiar example is $f(x)= e^{-x^2/2}$.    Indeed, defining as in \cite{titrim}, \cite{edw} the familiar Jacobi theta-function as 
$$\psi(x)= \sum_{n=1}^\infty e^{-n^2\pi x}, \  x >0,\eqno(3.5)$$
we know that it satisfies the functional equation 
$$\frac{1+2\psi(x)}{1+2 \psi(1/x)}= {1\over \sqrt x}.\eqno(3.6)$$ 
Now recalling Corollary 4 and equality (2.6), we obtain an interesting expansion
$${1\over e^{\pi x} -1} = \sum_{n: \  \mu(n)\neq 0}  \psi(nx),\  x >0,\eqno(3.7)$$
which can be proved by the interchange of summation due to the estimate
$$\sum_{n: \  \mu(n)\neq 0}  \sum_{m=1}^\infty e^{-n\pi m^2  x/2} e^{-n\pi m^2  x/2} \le  \sum_{m=1}^\infty e^{-m^2  x/2}   \sum_{n=1}^\infty e^{-n x/2},\ x >0 $$
and the use of (2.6).  Moreover, since \cite{titrim}
$$\psi(x)= {1\over 2\pi i} \int_{\sigma-i\infty}^{\sigma+ i\infty} \Gamma(s) \zeta(2s) (\pi x) ^{-s} ds,\  \sigma > 1,$$
the whole series 
$$\Psi(x) = (\Theta \psi)(x)= \sum_{n=1}^\infty   \psi(nx) = \sum_{n=1}^\infty  {1\over e^{n^2 \pi x} -1} ,\  x >0\eqno(3.8)$$ 
can be easily calculated by simple substitution and the interchange of summation and integration via the absolute convergence.  Hence  with (1.5) we get
$$\Psi(x)= {1\over 2\pi i} \int_{\sigma-i\infty}^{\sigma+ i\infty} \Gamma(s) \zeta(2s) \zeta(s) (\pi x) ^{-s} ds,\  \sigma > 1.\eqno(3.9)$$
Moving the contour to the left and taking into account the corresponding residues   in simple poles $s=1,  1/2$, we prove the following M\"{u}ntz type formulas 
$$\Gamma(s) \zeta(2s) \zeta(s) \pi^{-s} = \int_0^\infty x^{s-1} \left[ \Psi(x)- {\pi \over 6\  x} \right] dx,\   1/2 < \sigma < 1,
\eqno(3.10)$$
$$\Gamma(s) \zeta(2s) \zeta(s) \pi^{-s} = \int_0^\infty x^{s-1} \left[ \Psi(x)- {1\over \sqrt x} \  
\left( {\pi \over 6\ \sqrt x } + \zeta(1/2)\right) \   \right] dx,\  0  < \sigma < 1/2.\eqno(3.11)$$
Now, taking (3.8), we calculate the sum
$$\sum_{m=1}^\infty  \mu(m)  \Psi(mx),$$
 using the Lambert type expansion \cite{yakrim} of the exponential function, i.e.
 $$e^{-x}= \sum_{n=1}^\infty  {\mu(n) \over e^{n x} -1}.$$
 We have
 $$\sum_{m=1}^\infty  \mu(m)  \sum_{n=1}^\infty  {1\over e^{n^2 \pi m x} -1} =
 \sum_{n =1}^\infty  \sum_{m =1}^\infty  {\mu(m)  \over e^{n^2 \pi m x} -1} =
 \sum_{n =1}^\infty   e^{- n^2 \pi x} = \psi(x).
 $$
Thus, reciprocally,  (see (1.18), (3.8) ) 
$$\psi(x)=  (\Theta^{-1} \Psi)(x)= \sum_{m=1}^\infty  \mu(m)  \Psi(mx),\  x >0.\eqno(3.12)$$

Let us call $\hat{\psi} (x)$ the reduced  theta-function, which is defined by the series    (cf. (3.6))
$$\hat{\psi} (x)= \sum_{n: \mu(n)\neq 0}  e^{-n^2\pi x}, \  x >0. \eqno(3.13)$$
Since $e^{-x^2/2}$ is the self-reciprocal Fourier cosine transformation,  the corresponding formulas (2.14),  (2.15)  read in terms of the theta-functions (3.5), namely, 
 $$\hat{\psi} \left({x^2\over 2\pi}\right)  =   \sum_{n:\  \mu(n)\neq 0}    e^{-n^2 x^2/2}  =  \sum_{n =1}^\infty \mu(n) \psi \left({n^4x^2\over 2\pi}\right)$$
 $$= {3 \sqrt 2\over  x \pi\sqrt \pi}   
+  \sum_{n =1}^\infty \mu(n) \left[   {\sqrt{2\pi} \over n^2 x} \   \psi \left({2\pi\over n^4 x^2}\right)-   {1\over 2}\right], \   x >0. $$
Equivalently,  it can be written in the form (see (3.13))
 $$\hat{\psi} \left(x\right)  =    \sum_{n =1}^\infty \mu(n) \psi \left( n^4x\right)= {3 \over \pi^2 \sqrt  x }   
+  \sum_{n =1}^\infty \mu(n) \left[   {1 \over n^2 \sqrt x} \   \psi \left({1 \over n^4 x}\right)-   {1\over 2}\right], \   x >0.\eqno(3.14) $$
Moreover, as a consequence of the representation (2.3) and the Mellin transform formula (1.1), we have, correspondingly,
$$\pi^{-s/2} \Gamma\left({s\over 2}\right) {\zeta(s)\over \zeta(2s) } = \int_0^\infty x^{s/2-1}  \hat{\psi}(x) dx,\   \sigma > 1.\eqno(3.15)$$ 
Hence, following the Riemann technique (see, \cite{titrim}, Section 2.6), \cite{edw},  we use (3.14) to derive similarly the equalities 
$$\pi^{-s/2} \Gamma\left({s\over 2}\right) {\zeta(s)\over \zeta(2s) } =   \int_0^1  x^{s/2-1}  \hat{\psi}(x) dx +
 \int_1^\infty x^{s/2-1}  \hat{\psi}(x) dx= {6 \over \pi^2 (s-1)} $$
 $$  +  \int_0^1  x^{s/2-1}   \sum_{n =1}^\infty \mu(n) \left[   {1 \over n^2 \sqrt x} \   \psi \left({1 \over n^4 x}\right)-   {1\over 2}\right] dx +   \int_1^\infty x^{s/2-1} \sum_{n =1}^\infty \mu(n) \psi \left( n^4x\right)  dx$$ 
$$=  {6 \over \pi^2 (s-1)}   +  \int_1^\infty   x^{-s/2-1}   \sum_{n =1}^\infty \mu(n) \left[   {\sqrt x  \over n^2 } \   \psi \left({x\over n^4 }\right)-   {1\over 2}\right] dx +   \int_1^\infty x^{s/2-1} \sum_{n =1}^\infty \mu(n) \psi \left( n^4x\right)  dx.\eqno(3.16)$$ 
But the latter integral is an entire function.  It can be shown  via the estimate
$$\left| \int_1^\infty x^{s/2-1} \sum_{n =1}^\infty \mu(n) \psi \left( n^4x\right)  dx\right| \le 
\int_1^\infty x^{\sigma/2-1} \sum_{n =1}^\infty  \psi \left( n^4x\right)  dx $$
$$\le  \int_1^\infty x^{\sigma /2-1}  e^{-\pi x /3} dx \  \sum_{n =1}^\infty e^{-n^4 \pi/3}   \sum_{m=1}^\infty e^{-m^2 \pi/3} 
< \infty, \quad  \sigma \in \mathbb{R}.$$
Concerning another integral, we   first  recall the functional equation (3.6) to find  
$${1\over \sqrt x} \psi\left({1\over x}\right) - {1\over 2}= \psi(x)- {1\over 2\sqrt x}$$
and $\psi(x)=O\left( {1\over 2}( x^{-1/2}-1  )\right),\ x \to 0$. Hence,
$$ {6 \over \pi^2 (s-1)}   +  \int_1^\infty   x^{-s/2-1}   \sum_{n =1}^\infty \mu(n) \left[   {\sqrt x  \over n^2 } \   \psi \left({x\over n^4 }\right)-   {1\over 2}\right] dx =  {6 \over \pi^2 (s-1)}   +  \int_1^\infty   x^{-s/2-1}   \sum_{n =1}^\infty \mu(n) \left[   \psi \left({n^4\over x }\right)-   {\sqrt x\over 2 n^2}\right] dx $$
$$=  {6 \over \pi^2 (s-1)}   +  \int_1^\infty   x^{-s/2-1}   \sum_{n =1}^\infty \mu(n)  \psi \left({n^4\over x }\right) dx-
  {1\over 2}  \int_1^\infty   x^{-(s+1)/2}   dx  \sum_{n =1}^\infty {\mu(n)\over n^2} $$
$$=  \int_1^\infty   x^{-s/2-1}   \sum_{n =1}^\infty \mu(n)  \psi \left({n^4\over x }\right) dx,$$
where the splitting of the series is permitted due to them convergence. In fact, since 
$$\psi \left({n^4\over x }\right) = {1\over 4\pi i} \int_{\nu -i\infty}^{\nu+ i\infty}  \pi^{-w/2} \Gamma\left({w\over 2}\right) \zeta(w) n^{-2w} x^{w/2} dw,\  \nu > 1,$$
we have  by straightforward calculations 
$$ \int_1^\infty   x^{-s/2-1}   \sum_{n =1}^\infty \mu(n)  \psi \left({n^4\over x }\right) dx = 
{1\over 2\pi i} \int_{\nu -i\infty}^{\nu+ i\infty}  \frac{\pi^{-w/2} \Gamma\left({w\over 2}\right) \zeta(w)} {(s-w) \zeta(2w)} dw,
\ \sigma > \nu >1,$$
where   the interchange of summation and integration is possible owing to the absolute convergence.   Moving the contour in  the latter integral to the left and taking into account the pole of the integrand at $w=1$,  we establish finally the representation 
$$\pi^{-s/2} \Gamma\left({s\over 2}\right) {\zeta(s)\over \zeta(2s) } =   {6 \over \pi^2 (s-1)} + {1\over 2\pi i} \int_{\nu -i\infty}^{\nu+ i\infty}  \frac{\pi^{-w/2} \Gamma\left({w\over 2}\right) \zeta(w)} {(s-w) \zeta(2w)} dw$$
$$+  \int_1^\infty x^{s/2-1} \sum_{n =1}^\infty \mu(n) \psi \left( n^4x\right)  dx,\  {1\over 2} < \nu < 1.\eqno(3.17)$$
But the Cauchy type integral in (3.17) is analytic in any domain not containing points of the line $(\nu-i\infty,\ \nu +i\infty)$. Therefore equality (3.17) is valid for $s= 1/2 +it$.   Moreover,  via the Hadamard result $\zeta(1+2it)\neq 0,\  t \in \mathbb{R}$ and the   right -hand side of (3.17) is equal to zero when $s=1/2$.  Thus  we derive 
$$  {1\over 2\pi i} \int_{\nu -i\infty}^{\nu+ i\infty}  \frac{\pi^{-w/2} \Gamma\left({w\over 2}\right) \zeta(w)} {(1/2-w) \zeta(2w)} dw+  \int_1^\infty x^{-3/4} \sum_{n =1}^\infty \mu(n) \psi \left( n^4x\right)  dx =  {12\over \pi^2}.$$
Further, comparing  (3.17) with the Riemann equality \cite{titrim}
$$\pi^{-s/2} \Gamma\left({s\over 2}\right) \zeta(s) =  {1\over s(s-1)} + \int_1^\infty 
\left[ x^{s/2-1} + x^{-(s+1)/2} \right]  \psi (x) dx,$$ 
we observe,  that all zeros of its right-hand side, which  correspond to zeros of $\zeta(s)$ when $\sigma=1/2$ are zeros  on the critical line of the right-hand side of our equality (3.17), having one more zero at $s=1/2$.    By famous Hardy result about the infinity of zeros of $\zeta(s)$ on $\sigma=1/2$, we also conclude that the right-hand side of (3.17) has an infinity of zeros on the critical line.  Finally, putting $s=1/2+it$ in (3.17) one can conjecture that all zeros of its right -hand side are real,  giving  another  equivalence to   the Riemann hypothesis.

\bigskip
\centerline{{\bf Acknowledgments}}
\bigskip
The present investigation was supported, in part,  by the "Centro de
Matem{\'a}tica" of the University of Porto.

\bibliographystyle{amsplain}
{}

\end{document}